# HIDING A DRIFT

BY MIKLÓS RÁSONYI[1], WALTER SCHACHERMAYER[2] AND
RICHARD WARNUNG[2]

*Hungarian Academy of Sciences, University of Vienna and
Raiffeisen Capital Management*

In this article we consider a Brownian motion with drift of the form

$$dS_t = \mu_t\, dt + dB_t \qquad \text{for } t \geq 0,$$

with a specific nontrivial $(\mu_t)_{t\geq 0}$, predictable with respect to $\mathbb{F}^B$, the natural filtration of the Brownian motion $B = (B_t)_{t\geq 0}$. We construct a process $H = (H_t)_{t\geq 0}$, also predictable with respect to $\mathbb{F}^B$, such that $((H \cdot S)_t)_{t\geq 0}$ is a Brownian motion in its own filtration. Furthermore, for any $\delta > 0$, we refine this construction such that the drift $(\mu_t)_{t\geq 0}$ only takes values in $]\mu - \delta, \mu + \delta[$, for fixed $\mu > 0$.

**1. Introduction.** Let $B = (B_t)_{t\geq 0}$ be a standard Brownian motion on a probability space $(\Omega, \mathcal{F}, \mathbb{P})$. For a fixed constant $\mu > 0$, denote the Brownian motion with drift $\mu$ as $S = (S_t)_{t\geq 0}$, defined by

(1.1) $$S_t = \mu t + B_t, \qquad t \geq 0.$$

Furthermore, let $\mathbb{F}^B := (\mathcal{F}^B_t)_{t\geq 0}$ denote the right continuous, saturated filtration generated by $B$. Given a predictable, $\mathbb{F}^B$-adapted process $H = (H_t)_{t\geq 0}$, we consider the stochastic integral $(H \cdot S) = ((H \cdot S)_t)_{t\geq 0}$ in its right continuous, saturated filtration $\mathbb{F}^{(H\cdot S)} := (\mathcal{F}^{(H\cdot S)}_t)_{t\geq 0}$. Marc Yor posed the following question:

Received February 2008; revised January 2009.
[1]Supported by Austrian Science Fund (FWF) Grant P19456 and from the Hungarian Science Foundation (OTKA) Grant F049094.
[2]Supported by Austrian Science Fund (FWF) Grant P19456, Vienna Science and Technology Fund (WWTF) Grant MA13 and the Christian Doppler Research Association (CDG). While the research of this paper was conducted all three authors were affiliated with the Institute for Mathematical Methods in Economics, Vienna University of Technology.

**Question 1**: Can we define an $\mathbb{F}^B$-predictable process $H$ such that the resulting stochastic integral $(H \cdot S)$ is a Brownian motion (without drift) in its own filtration, that is, an $\mathbb{F}^{(H \cdot S)}$-Brownian motion?

Clearly, the predictable integrand $H$ can only take values in $\{-1, 1\}, \mathbb{P} \otimes \lambda$-a.s., $\lambda$ denoting Lebesgue measure on $[0, \infty)$, in order to make sure that the process $(H \cdot S)$ has the quadratic variation of a Brownian motion.

In fact, at first glance it seems completely unlikely that an $\mathbb{F}^B$-predictable process $H$ with the required property does exist. Indeed, intuitively speaking, it would have to start with $\mathbb{P}[H_0 = 1] = \mathbb{P}[H_0 = -1] = 1/2$, which seems absurd, as $H_0$ is required to be $\mathcal{F}_0^B$-measurable and therefore $\mathbb{P}$-a.s. constant (the sigma-algebra $\mathcal{F}_0^B$ is trivial). Fortunately this intuitive argument is not quite correct, as the predictable process $H = (H_t)_{t \geq 0}$ is only defined modulo $\mathbb{P} \otimes \lambda$ null-sets, so that it does not really makes sense to speak about the random variable $H_0$. Nevertheless, the preceding heuristics seem to indicate that we need some random sign $\varepsilon$ with $\mathbb{P}[\varepsilon = 1] = \mathbb{P}[\varepsilon = -1] = 1/2$ which is *independent* of the Brownian motion $B$ to be able to start a successful construction of the desired integrand $H = (H_t)_{t \geq 0}$ for $t$ close to $t = 0$.

So, let us cheat for a moment and fix a random variable $\varepsilon$, defined on $(\Omega, \mathcal{F}, \mathbb{P})$ with $\mathbb{P}[\varepsilon = 1] = \mathbb{P}[\varepsilon = -1] = 1/2$, and consider the enlarged filtration $\mathbb{F}^{B, \varepsilon}$ defined by letting $\mathcal{F}_t^{B, \varepsilon} = \sigma(\mathcal{F}_t^B, \varepsilon)$ for $t \geq 0$.

Let us now try to construct an integrand $H = (H_t)_{t \geq 0}$ which is predictable in the enlarged filtration $\mathbb{F}^{B, \varepsilon}$ and such that the stochastic integral $(H \cdot S)$ is a Brownian motion (without drift) in its own filtration $\mathbb{F}^{(H \cdot S)}$. We have an obvious way to start the construction of $H$ at time $t = 0$ by letting

(1.2) $$H_0 := \varepsilon,$$

or rather, reasoning heuristically with infinitesimals,

$$H_u := \varepsilon \qquad \text{for } 0 \leq u \leq dt.$$

This yields an integrand $(H_u)_{0 \leq u \leq dt}$ such that the stochastic integral $(H \cdot S)_{0 \leq u \leq dt}$ is a martingale for the infinitesimal time interval $[0, dt]$. Indeed,

$$\mathbb{E}[d(H \cdot S)_0] = \mathbb{E}[\varepsilon(S_{dt} - S_0)] = \mathbb{E}[\varepsilon(B_{dt} - B_0) + \varepsilon \mu \, dt] = 0 \, dt.$$

But already an infinitesimal instant of time later we again are in trouble: after having observed the process $(H \cdot S)$ during the infinitesimal time interval $[0, dt]$, we have learned something (which turns out to be of the order $dt^{1/2}$) on the probability of $\varepsilon$ equaling $+1$ or $-1$, conditionally on the process $(H \cdot S)_{0 \leq u \leq dt}$.

Hence, the approach of defining $H_t = H_0 = \varepsilon$ for $t \in [0, \Delta t]$ for a *finite* increment $\Delta t > 0$ yields a process $(H \cdot S)_{0 \leq t \leq \Delta t}$ which *fails to be* a martingale in its own filtration, as one easily verifies.

At this stage we remembered Pólya's famous dictum:

"To every problem there is an easier problem."

Instead of asking Yor's original question for the process $S$ with constant drift $\mu$, as in (1.1) we pose the same question, but with $\mu$ replaced by an appropriate predictable process $(\mu_t)_{t \geq 0}$, that is,

$$dS_t = \mu_t \, dt + dB_t, \tag{1.3}$$

where $(\mu_t)_{t \geq 0}$ is tailor-made such that, for the integrand $H_t = \varepsilon$, for $t \geq 0$, we indeed obtain a process $(H \cdot S)_{t \geq 0}$ which is a Brownian motion in its own filtration $\mathbb{F}^{(H \cdot S)}$. This program indeed turns out to be doable as summarized in the subsequent statement.

PROPOSITION 1.1. *Suppose that on $(\Omega, \mathcal{F}, \mathbb{P})$ there is a standard Brownian motion $B = (B_t)_{t \geq 0}$ and a random variable $\varepsilon$ with $\mathbb{P}[\varepsilon = 1] = \mathbb{P}[\varepsilon = -1] = 1/2$, independent of $B$. Denote by $\mathbb{F}^B$ the filtration generated by $B$.*

*For each $\mu > 0$, there is an $\mathbb{F}^B$-predictable process $\mu_t$ taking values in $]0, 2\mu[$ such that defining $S = (S_t)_{t \geq 0}$ by $S_0 = 0$ and*

$$dS_t = \mu_t \, dt + dB_t, \qquad t \geq 0, \tag{1.4}$$

*we have that*

$$Y_t = \varepsilon S_t, \qquad t \geq 0,$$

*is a Brownian motion in its own filtration.*

The preceding result is a preliminary step toward a satisfactory answer to Yor's question. It has two deficiencies: first, we had to replace the constant $\mu$ by a process $(\mu_t)_{t \geq 0}$ fluctuating in $]0, 2\mu[$, and, second, we had to enlarge the filtration $\mathbb{F}^B$ to $\mathbb{F}^{\varepsilon,B}$ in order to be able to define our predictable integrand $H_t \equiv \varepsilon$, for $t \geq 0$.

As regards the second issue, we can get completely rid of the necessity of introducing the additional source of randomness $\varepsilon$ by applying the Lévy transform; see Section 3. We can indeed find an integrand $H$ which is predictable with respect to $\mathbb{F}^B$ instead of $\mathbb{F}^{\varepsilon,B}$ and still does the job. As regards the first issue, we can refine the construction in such a way that the process $(\mu_t)_{t \geq 0}$ only takes values in $]\mu - \delta, \mu + \delta[$ instead of $]0, 2\mu[$, for given $\delta > 0$.

We summarize our findings in the subsequent theorem. Very roughly speaking, it states that the answer to Yor's question is positive, provided that we allow for a tailor-made drift process $(\mu_t)_{t \geq 0}$ instead of a constant drift $\mu$, which may be chosen to satisfy $|\mu_t - \mu| < \delta$, for given $\delta > 0$.

THEOREM 1.2. *Let $B = (B_t)_{t \geq 0}$ be a Brownian motion defined on the filtered probability space $(\Omega, \mathcal{F}, \mathbb{F}^B, \mathbb{P})$, where $\mathbb{F}^B = (\mathcal{F}^B_t)_{t \geq 0}$ is the right-continuous saturated filtration generated by $B$:*

(i) *For each $\mu > 0$, there are $\mathbb{F}^B$-predictable processes $H = (H_t)_{t \geq 0}$, taking values in $\{-1, 1\}$, and $(\mu_t)_{t \geq 0}$, taking values in $]0, 2\mu[$, such that for $S = (S_t)_{t \geq 0}$ defined by $S_0 = 0$ and*

$$dS_t = \mu_t \, dt + dB_t$$

*we have that the process $((H \cdot S)_t)_{t \geq 0}$ is a Brownian motion in its own filtration $\mathbb{F}^{(H \cdot S)}$.*

(ii) *Furthermore, for each $\delta > 0$ we can choose $(\mu_t)_{t \geq 0}$ such that it only takes values in $]\mu - \delta, \mu + \delta[$.*

However, Question 1 of M. Yor in its original form above still remains an open and challenging problem. For recent work on the conservation of the martingale property under a change of filtration, we refer to [1].

M. Émery asked us the following question: what about the discrete-time version of the problem? The proper discrete analogue is an i.i.d. sequence $(\varepsilon_n)_{n \leq 0}$ in its natural filtration $(\mathcal{F}_n)_{n \leq 0}$ such that $\mathbb{P}[\varepsilon_n = 1] = 1 - \mathbb{P}[\varepsilon_n = -1] = p \in ]0, 1[ \setminus \{1/2\}$. The question now reads as follows: is there an $(\mathcal{F}_n)_{n \leq 0}$-predictable sequence $(h_n)_{n \leq 0}$ of $\{-1, 1\}$-valued random variables such that the sequence $(h_n \varepsilon_n)_{n \leq 0}$ is i.i.d. with $\mathbb{P}[h_n \varepsilon_n = 1] = \mathbb{P}[h_n \varepsilon_n = -1] = 1/2$?

This discrete version turns out to be simpler than the continuous one and we shall give in the Appendix a positive solution, even in a slightly more general setting.

The article is structured as follows: In Section 2 we describe the details of the construction of $(\mu_t)_{t \geq 0}$ with respect to $\mathbb{F}^{\varepsilon, B}$. Next, in Section 3 we prove the first part of Theorem 1.2 above. Finally, in Section 4 we use stopping techniques in order to show the second statement of Theorem 1.2.

**2. Constructing the drift process.** Fix a probability space $(\Omega, \mathcal{F}, P)$. Let $B$ be a Brownian motion and let $\varepsilon$ be an independent random sign with $P(\varepsilon = 1) = P(\varepsilon = -1) = 1/2$. Consider

$$(2.1) \qquad S_t := \int_0^t \mu_s \, ds + B_t, \qquad t \geq 0,$$

with some bounded $\mathbb{F}^B$-predictable drift $\mu_t$ and set $Y_t := \varepsilon S_t$. Our purpose is to find $\mu_t$ such that $Y_t$ is a Brownian motion in its own filtration.

We imagine $\mu_t$ as being "glued together" from two $\mathbb{F}^Y$-predictable processes. Formally, let $\mu_t^+, \mu_t^-$ be $\mathbb{F}^Y$-predictable bounded processes such that

$$(2.2) \qquad \mu_t := 1_{\{\varepsilon = +1\}} \mu_t^+ + 1_{\{\varepsilon = -1\}} \mu_t^-.$$

We wish to derive conditions on $\mu_t^+, \mu_t^-$ which ensure that $Y_t$ is as required. To this end, introduce the conditional probabilities

$$(2.3) \qquad p_t := P[\varepsilon = 1 | \mathcal{F}_t^Y], \qquad t \geq 0.$$

In the language of filtering theory $p_t$ gives the distribution of the "signal" $\varepsilon$, conditionally on the "observations" $Y$.

PROPOSITION 2.1. *Let $S$ and $Y$ be as above and $(\mu_t)_{t\geq 0}$ a bounded $\mathbb{F}^B$-predictable process of the form (2.2). The conditional probabilities $p_t$ defined in (2.3) satisfy $p_0 = 1/2$ and*

$$dp_t = \frac{(\mu_t^+ + \mu_t^-)^2}{2}[\varepsilon p_t(1-p_t) + p_t(1-p_t)(1-2p_t)]\,dt \quad (2.4)$$
$$+ \varepsilon p_t(1-p_t)(\mu_t^+ + \mu_t^-)\,dB_t.$$

PROOF. For each $T > 0$ there exists a measure $\mathbb{Q}_T \sim \mathbb{P}|_{\mathcal{F}_T^{\varepsilon,B}}$ such that $(Y_t)_{0 \leq t \leq T}$ is a $(\mathbb{Q}_T, \mathbb{F}^{\varepsilon,B})$-Brownian motion. By the Girsanov theorem, we know that the Radon–Nikodym derivative is given by

$$(2.5) \quad \frac{d\mathbb{Q}_T}{d\mathbb{P}} = \begin{cases} \exp\left(-\int_0^T \mu_t^+ \, dB_t - (1/2)\int_0^T (\mu_t^+)^2\,dt\right), & \text{if } \varepsilon = 1, \\ \exp\left(-\int_0^T \mu_t^- \, dB_t - (1/2)\int_0^T (\mu_t^-)^2\,dt\right), & \text{if } \varepsilon = -1. \end{cases}$$

It follows that, for each $T > 0$, the $(\mathbb{Q}_T, \mathbb{F}^{\varepsilon,B})$-martingale $Z_t := \frac{d\mathbb{P}}{d\mathbb{Q}_T}|_{\mathcal{F}_t^{\varepsilon,B}}, 0 \leq t \leq T$, is of the form

$$(2.6) \qquad Z_t = Z_t^+ 1_{\{\varepsilon=1\}} + Z_t^- 1_{\{\varepsilon=-1\}}, \qquad 0 \leq t \leq T,$$

where the processes $(Z_t^+)_{0 \leq t \leq T}$ and $(Z_t^-)_{0 \leq t \leq T}$ are given by $Z_0^+ = Z_0^- = 1$ and

$$dZ_t^+ = \mu_t^+ Z_t^+ \, dY_t,$$
$$dZ_t^- = -\mu_t^- Z_t^- \, dY_t,$$

respectively. Note that $\mu_t^+, \mu_t^-$ are assumed to be $\mathbb{F}^Y$-predictable and, thus, $(Z_t^+)_{0 \leq t \leq T}$ and $(Z_t^-)_{0 \leq t \leq T}$ are $\mathbb{F}^Y$-predictable, too. By the assumption on $\mu_t$, $Z_t$ is clearly $\mathbb{F}^B$-predictable.

We claim that, under $\mathbb{Q}_T$, $\varepsilon$ is independent of $Y$ and has the same law as under $\mathbb{P}$. Indeed, as $Z_T$ is $\mathcal{F}_T^B$-measurable and $\varepsilon$ is independent of $B$ (and hence of $S$) under $\mathbb{P}$, for any bounded measurable functions $h, j$ we have

$$\mathbb{E}_{\mathbb{Q}_T}[h(S)j(\varepsilon)] = \mathbb{E}_{\mathbb{P}}[(1/Z_T)h(S)]\mathbb{E}_{\mathbb{P}}[j(\varepsilon)] = \mathbb{E}_{\mathbb{Q}_T}[h(S)]\mathbb{E}_{\mathbb{Q}_T}[j(\varepsilon)],$$

showing the $\mathbb{Q}_T$-independence of $S$ and $\varepsilon$ as well as $\mathbb{Q}_T[\varepsilon = \pm 1] = 1/2$.

Now note that, under $\mathbb{Q}_T$, $S$ is a Brownian motion. By symmetry of the Brownian motion, we also get that $Y = \varepsilon S$ and $\varepsilon$ are $\mathbb{Q}_T$-independent as claimed above.

Clearly, we have $1_{\{\varepsilon=1\}} = \frac{\varepsilon+1}{2}$. Thus, for calculating

$$\mathbb{E}_{\mathbb{P}}[1_{\{\varepsilon=1\}} \mid \mathcal{F}_t^Y] = \tfrac{1}{2}(\mathbb{E}_{\mathbb{P}}[\varepsilon \mid \mathcal{F}_t^Y] + 1), \qquad t \geq 0,$$

we first calculate $\mathbb{E}_{\mathbb{P}}[\varepsilon \mid \mathcal{F}_t^Y], t \geq 0$. Fix any $T > 0$ and consider that by Bayes' formula and the tower law applied to the $(\mathbb{Q}_T, \mathbb{F}^{\varepsilon,B})$-martingale $Z_T$, it holds that

$$\mathbb{E}_{\mathbb{P}}[\varepsilon \mid \mathcal{F}_t^Y] = \frac{\mathbb{E}_{\mathbb{Q}_T}[\varepsilon Z_T \mid \mathcal{F}_t^Y]}{\mathbb{E}_{\mathbb{Q}_T}[Z_T \mid \mathcal{F}_t^Y]} = \frac{\mathbb{E}_{\mathbb{Q}_T}[\varepsilon \mathbb{E}_{\mathbb{Q}_T}[Z_T \mid \mathcal{F}_t^{\varepsilon,B}] \mid \mathcal{F}_t^Y]}{\mathbb{E}_{\mathbb{Q}_T}[\mathbb{E}_{\mathbb{Q}_T}[Z_T \mid \mathcal{F}_t^{\varepsilon,B}] \mid \mathcal{F}_t^Y]}$$
$$= \frac{\mathbb{E}_{\mathbb{Q}_T}[\varepsilon Z_t \mid \mathcal{F}_t^Y]}{\mathbb{E}_{\mathbb{Q}_T}[Z_t \mid \mathcal{F}_t^Y]}, \qquad 0 \leq t \leq T.$$

By independence of $\varepsilon$ and $Y$ (under $\mathbb{Q}_T$) and $\mathbb{F}^Y$-adaptedness of $(Z_t^+)_{0 \leq t \leq T}$ and $(Z_t^-)_{0 \leq t \leq T}$, we get

$$\frac{\mathbb{E}_{\mathbb{Q}_T}[\varepsilon Z_t \mid \mathcal{F}_t^Y]}{\mathbb{E}_{\mathbb{Q}_T}[Z_t \mid \mathcal{F}_t^Y]} = \frac{\mathbb{E}_{\mathbb{Q}_T}[Z_t^+ 1_{\{\varepsilon=1\}} \mid \mathcal{F}_t^Y] - \mathbb{E}_{\mathbb{Q}_T}[Z_t^- 1_{\{\varepsilon=-1\}} \mid \mathcal{F}_t^Y]}{\mathbb{E}_{\mathbb{Q}_T}[Z_t^+ 1_{\{\varepsilon=1\}} \mid \mathcal{F}_t^Y] + \mathbb{E}_{\mathbb{Q}_T}[Z_t^- 1_{\{\varepsilon=-1\}} \mid \mathcal{F}_t^Y]} = \frac{Z_t^+ - Z_t^-}{Z_t^+ + Z_t^-}$$
$$= \left[\exp\left(\int_0^t \mu_u^+ \, dY_u - \frac{1}{2}\int_0^t (\mu_u^+)^2 \, du\right)\right.$$
$$\left. - \exp\left(-\int_0^t \mu_u^- \, dY_u - \frac{1}{2}\int_0^t (\mu_u^-)^2 \, du\right)\right]$$
$$\times \left[\exp\left(\int_0^t \mu_u^+ \, dY_u - \frac{1}{2}\int_0^t (\mu_u^+)^2 \, du\right)\right.$$
$$\left. + \exp\left(-\int_0^t \mu_u^- \, dY_u - \frac{1}{2}\int_0^t (\mu_u^-)^2 \, du\right)\right]^{-1}$$
$$= \frac{\exp(\int_0^t (\mu_u^+ + \mu_u^-) \, dY_u - (1/2)\int_0^t [(\mu_u^+)^2 - (\mu_u^-)^2] \, du) - 1}{\exp(\int_0^t (\mu_u^+ + \mu_u^-) \, dY_u - (1/2)\int_0^t [(\mu_u^+)^2 - (\mu_u^-)^2] \, du) + 1},$$
$$0 \leq t \leq T.$$

So we have

$$\mathbb{E}_{\mathbb{P}}[1_{\{\varepsilon=1\}} \mid \mathcal{F}_t^Y]$$
(2.7)
$$= \frac{1}{2}\left(\frac{\mathbb{E}_{\mathbb{Q}_T}[\varepsilon Z_t \mid \mathcal{F}_t^Y]}{\mathbb{E}_{\mathbb{Q}_T}[Z_t \mid \mathcal{F}_t^Y]} + 1\right)$$
$$= \frac{\exp(\int_0^t (\mu_u^+ + \mu_u^-) \, dY_u - (1/2)\int_0^t [(\mu_u^+)^2 - (\mu_u^-)^2] \, du)}{\exp(\int_0^t (\mu_u^+ + \mu_u^-) \, dY_u - (1/2)\int_0^t [(\mu_u^+)^2 - (\mu_u^-)^2] \, du) + 1},$$

for $0 \leq t \leq T$. Define the process $(U_t)_{t \geq 0}$ given by $U_0 = 0$ and

$$dU_t = (\mu_u^+ + \mu_u^-) \, dY_t - \frac{(\mu_u^+)^2 - (\mu_u^-)^2}{2} \, dt.$$

Applying the Itô formula to (2.7) and recalling the expression for $(Y_t)_{t\geq 0}$, we get

$$d\frac{\exp(U_t)}{\exp(U_t)+1} = \frac{\exp(U_t)}{(\exp(U_t)+1)^2}\left(\varepsilon\frac{(\mu_t^+ + \mu_t^-)^2}{2}\,dt + \varepsilon(\mu_t^+ + \mu_t^-)\,dB_t\right)$$
$$+ \frac{1}{2}\frac{\exp(U_t)-\exp(2U_t)}{(\exp(U_t)+1)^3}(\mu_t^+ + \mu_t^-)^2\,dt.$$

Using $p_t = \frac{\exp(U_t)}{\exp(U_t)+1}, t\geq 0$, we get (2.4). □

We also have the following proposition.

PROPOSITION 2.2. *Under the assumptions of Proposition 2.1, suppose, in addition, that for all $u \geq 0$,*

(2.8) $$p_u\mu_u^+ - (1-p_u)\mu_u^- = 0 \quad \text{a.s.,}$$

*then the process $Y$ is an $\mathbb{F}^Y$-Brownian motion.*

PROOF. Obviously $(Y_t)_{t\geq 0}$ is $\mathbb{F}^Y$-adapted, continuous and has the right quadratic variation as the drift is bounded. In order to fulfill Lévy's characterization theorem of Brownian motion, we need to check the martingale condition. Therefore, fix $s \leq t < \infty$ and consider that

$$\mathbb{E}[\varepsilon(S_t - S_s) \mid \mathcal{F}_s^Y] = \mathbb{E}\left[\int_s^t \varepsilon\mu_u\,du \,\Big|\, \mathcal{F}_s^Y\right] + \mathbb{E}[\varepsilon(B_t - B_s) \mid \mathcal{F}_s^Y].$$

The second conditional expectation is 0. The martingale property is thus equivalent to

(2.9) $$\mathbb{E}\left[\int_s^t \varepsilon\mu_u\,du\, 1_A\right] = 0,$$

for all $A \in \mathcal{F}_s^Y$. Note that the Fubini theorem applies as $|\varepsilon\mu_u 1_A|$ is bounded. Furthermore, using the tower law, we get that (2.9) holds iff

$$\mathbb{E}\left[\int_s^t \varepsilon\mu_u\,du\, 1_A\right] = \int_s^t \mathbb{E}[\mathbb{E}[\varepsilon\mu_u \mid \mathcal{F}_u^Y]1_A]\,du = 0,$$

for all $A \in \mathcal{F}_u^Y$. Recall that $\mu_u^+$ and $\mu_u^-$ are assumed to be $\mathcal{F}_u^Y$-measurable for $u \geq 0$. It follows from the hypotheses of this proposition that

$$\mathbb{E}[\varepsilon\mu_u \mid \mathcal{F}_u^Y] = p_u\mu_u^+ - (1-p_u)\mu_u^- = 0,$$

concluding the proof. □

Formula (2.8) shows that it is reasonable to choose $\mu_t^+$ proportional to $(1-p_t)$ and $\mu_t^-$ proportional to $p_t$. This will guarantee the validity of (2.8), as the next proposition shows.

PROPOSITION 2.3. *Let $\mu > 0$ be an arbitrary constant. Let $g_t$ be a solution of the equation*

$$(2.10) \quad dg_t = 2\mu^2[\varepsilon g_t(1-g_t) + g_t(1-g_t)(1-2g_t)] \, dt + 2\mu\varepsilon g_t(1-g_t) \, dB_t,$$

$$g_0 = 1/2,$$

*adapted to the filtration $\mathbb{F}^{\varepsilon,B}$ and satisfying $0 \leq g_t \leq 1$ for $t \geq 0$. Set*

$$(2.11) \quad \mu_t^+ = 2\mu(1-g_t), \qquad \mu_t^- = 2\mu g_t, \qquad t \geq 0,$$

*and define $\mu_t, S_t, Y_t, p_t$ accordingly. If $g_t$ is $\mathbb{F}^Y$-predictable and $\mu_t$ is $\mathbb{F}^B$-predictable, then $p_t$ equals $g_t$ and $Y$ is a Brownian motion in its own filtration.*

PROOF. First note that the coefficients of the autonomous SDE (2.10) are Lipschitz-continuous and bounded when restricted to the interval $[0,1]$, hence, $g_t$ is the *unique strong* solution of (2.10) adapted to $\mathbb{F}^{\varepsilon,B}$. Thus, it suffices to prove that $p_t$ is also a solution of (2.10). Obviously, $p_0 = \mathbb{P}[\varepsilon = 1] = 1/2 = g_0$.

If $g_t$ is $\mathbb{F}^Y$-predictable, then so are $\mu_t^+, \mu_t^-$. Proposition 2.1 shows that $p_t$ is a solution of (2.10), hence, indeed, $p_t = g_t$. With this choice of $\mu_t^+, \mu_t^-$, equation (2.8) is satisfied and Proposition 2.2 allows us to conclude. □

It remains to solve the stochastic differential equation (2.10).

PROOF OF PROPOSITION 1.1. In Section 4.2 of [4] a different filtering problem leads to almost the same equation as (2.10). That equation has an explicit solution [see (4.55) on p. 180] from which it is easy to make the guess

$$(2.12) \quad g_t := \begin{cases} \dfrac{\exp(2\mu B_t + 2\mu^2 t)}{1 + \exp(2\mu B_t + 2\mu^2 t)}, & \text{if } \varepsilon = 1, \\ \dfrac{1}{1 + \exp(2\mu B_t + 2\mu^2 t)}, & \text{if } \varepsilon = -1. \end{cases}$$

Applying Itô's formula, we may check that this indeed gives a (strong) solution to (2.10) which trivially stays in $(0,1)$. Define $\mu_t^+, \mu_t^-, \mu_t, S_t, Y_t, p_t$ as in Proposition 2.3. One may check that

$$(2.13) \quad dg_t = 2\mu g_t(1-g_t) \, dY_t,$$

showing that $g_t$ is $\mathbb{F}^Y$-predictable. We find that the dynamics of $\mu_t$ is

$$(2.14) \quad d\mu_t = -\mu_t^2(2\mu - \mu_t) \, dt - \mu_t(2\mu - \mu_t) \, dB_t,$$

hence, $\mu_t$ is $\mathbb{F}^B$-predictable. For later use we note that, substituting in to (2.11), we get the following formula for $\mu_t$:

$$\mu_t = \frac{2\mu}{1 + \exp(2\mu B_t + 2\mu^2 t)}. \tag{2.15}$$

Proposition 2.3 now implies that $p_t = g_t$ and $Y_t$ is indeed as required. $\square$

**3. Passing to the Lévy transform.** In this section we describe how to get rid of the enlargement of the filtration $\mathbb{F}^B$ by the sign $\varepsilon$. We will make use of the Lévy transform which arises naturally in the famous Tanaka formula for the SDE of $(|B_t|)_{t \geq 0}$ for some Brownian motion $B$ (for the derivation of the Tanaka formula, see, e.g., [3]).

Recall that the Lévy transform $(M_t^0)_{t \geq 0}$ of a Brownian motion $(B_t)_{t \geq 0}$ is defined by

$$M_t^0 = \int_0^t \operatorname{sign}(B_s)\, dB_s, \qquad t \geq 0, \tag{3.1}$$

where we use the sign function in the following left-continuous form:

$$\operatorname{sign}(x) = 1_{\{x > 0\}} - 1_{\{x \leq 0\}} \qquad \text{for } x \in \mathbb{R}.$$

Among the properties of the Lévy transform, we mention that $(M_t^0)_{t \geq 0}$ is a Brownian motion in its own filtration and that the filtration generated by $(M_t^0)_{t \geq 0}$ equals the one generated by $(|B_t|)_{t \geq 0}$ which is strictly smaller than the filtration generated by $(B_t)_{t \geq 0}$.

PROOF OF (i) IN THEOREM 1.2. We come back to the setting of Section 2 and consider the filtered probability space $(\Omega, \mathcal{F}, \mathbb{F}^{\varepsilon, B}, \mathbb{P})$. Let us take $\mu_t, S_t$ as constructed in the proof of Proposition 1.1. Introduce $(Y_t)_{t \geq 0} = (\varepsilon S_t)_{t \geq 0}$; this is a Brownian motion in its own filtration, by Proposition 1.1.

Now consider the Lévy transform $(M_t)_{t \geq 0}$ of the $\mathbb{F}^Y$-Brownian motion $Y = (Y_t)_{t \geq 0}$. It is defined by $M_0 = 0$ and

$$\begin{aligned} dM_t &= \operatorname{sign}(Y_t)\, dY_t = \operatorname{sign}(\varepsilon S_t)\varepsilon\, dS_t \\ &= \operatorname{sign}(S_t)\, dS_t, \qquad t \geq 0. \end{aligned} \tag{3.2}$$

This is again a Brownian motion (in $\mathbb{F}^Y$ as well as in $\mathbb{F}^M$) by the properties of the Lévy-transform. It follows that, with the choice $H_t = \operatorname{sign}(S_t)$, the process $(H \cdot S)_t, t \geq 0$ is a Brownian motion in its own filtration. $\square$

**4. $L^\infty$-approximation of a constant drift.** The aim of this section is to show that we can in fact define a process $(S_t)_{t\geq 0}$ such that the drift is close to a constant drift $\mu$ with respect to the norm in $L^\infty$. The strategy is that we stop whenever the drift $(\mu_t)_{t\geq 0}$ has moved by some small fixed number. After that we will restart the construction. A somewhat delicate point is that the stopping has to be done in a way adapted to $\mathbb{F}^M$. Lemma 4.1 below shows that this is indeed possible.

The distance of the drift process $\mu_t$ from $\mu$ is proportional to the distance of $p_t$ from one half. Namely, by (2.2) and (2.11),

$$
\begin{aligned}
|\mu_t - \mu| &= |2\mu(1-p_t) - \mu| 1_{\{\varepsilon=1\}} + |2\mu p_t - \mu| 1_{\{\varepsilon=-1\}} \\
&= 2\mu|\tfrac{1}{2} - p_t| \qquad \text{for } t \geq 0.
\end{aligned}
\tag{4.1}
$$

In the following lemma we show that one can define a stopping time in the filtration generated by the Lévy transform $(M_t)_{t\geq 0}$ of the Brownian motion $(Y_t)_{t\geq 0}$, that is, $\mathbb{F}^M := (\mathcal{F}_t^M)_{t\geq 0}$, such that we have a control over the distance of $p$ from $1/2$.

LEMMA 4.1. *Take $p_t, \mu_t, S_t, Y_t$ as constructed in the proof of Proposition 1.1 in Section 2 and consider the Lévy transform $M = (M_t)_{t\geq 0}$ of $Y = (Y_t)_{t\geq 0}$, defined by $M_0 = 0$ and*

$$dM_t = \text{sign}(Y_t)\, dY_t.$$

*For each $\delta > 0$ we define the stopping time $\rho_\delta := \inf\{t : |M_t| \geq \delta\} \wedge \delta$. The following estimate holds:*

$$
\left| p_t - \frac{1}{2} \right| \leq \delta \frac{(2\mu + 3\mu^2)}{2} \qquad \text{for } 0 \leq t \leq \rho_\delta.
\tag{4.2}
$$

PROOF. We present the proof in two steps.

*Step 1*: We show that $|Y_t| \leq 2\delta$, for $0 \leq t \leq \rho_\delta$.

Let $\tau := \inf\{t : |Y_t| \geq 2\delta\}$ and let $\sigma := \max\{t \leq \tau : Y_t = 0\}$, that is, the time of the last zero of $Y$ preceding $\tau$. We note in passing that $\tau$ is a stopping time in the filtration $\mathbb{F}^Y$, while $\sigma$ fails to be a stopping time. Observe that by Tanaka's formula (see, for instance, [3])

$$M_t = |Y_t| - L_t,$$

where $L$ is the local time of $Y$ at zero. By definition of $\sigma$, the local time $L$ does not grow on $[\sigma, \tau]$ and, thus, a.s. $L_\sigma = L_\tau$. For the process $M$ this gives

$$|M_\sigma - M_\tau| = 2\delta,$$

so that $\sup_{0 \leq t \leq \tau} |M_t| \geq \delta$, which shows that a.s. $\rho_\delta \leq \tau$, that is, $|Y_t| \leq 2\delta$ for $0 \leq t \leq \rho_\delta$.

*Step 2*: By straightforward calculation, $|p_t - (1/2)| = (1/2)|\operatorname{th}(\mu^2 t + \mu B_t)|$, where th denotes the hyperbolic tangent. As $|\operatorname{th} x| \leq |x|$, $dY_t = \varepsilon \mu_t \, dt + \varepsilon \, dB_t$ and $\mu_t \in \,]0, 2\mu[$, we may write, for $t \leq \rho_\delta$,

$$\left| p_t - \frac{1}{2} \right| \leq \mu^2 \frac{t}{2} + \mu \frac{|B_t|}{2} \leq \mu^2 \frac{t}{2} + \mu \frac{|Y_t|}{2} + \mu \left| \int_0^t \mu_s \, ds \right| / 2 \leq \mu^2 \frac{\delta}{2} + \mu \frac{2\delta}{2} + \mu^2 \delta,$$

using Step 1 and $\rho_\delta \leq \delta$. $\square$

Using the previous lemma, we can refine the construction by stopping and restarting, when we are too far away from a constant drift, considering the information of $\mathbb{F}^M$ only. Fix the constant $\mu > 0$. For the goal of controlling the $L^\infty$ distance of $\mu$ and the drift process $\mu_t$ to be constructed, fix also a constant $\delta > 0$.

The strategy is straightforward. We start at $t = 0$ using the drift $(\mu_t^1)_{t \geq 0}$ which is given by $\mu_0^1 = \mu$ and (2.14). Define the process $S^1$ by $S_0^1 = 0$,

$$dS_t^1 = \mu_t^1 \, dt + dB_t, \qquad t \geq 0.$$

We perform the Lévy transform which results in a process $(M_t^1)_{t \geq 0}$. Introducing the $\mathbb{F}^{M^1}$-stopping time

$$\tau_1 := \inf\{t > 0 : |M_t^1| \geq \delta\} \wedge \delta,$$

we can assure by Lemma 4.1 and (4.1) that

(4.3) $\qquad |\mu_t^1 - \mu| \leq \delta(3\mu^3 + 2\mu^2) \qquad \text{for } 0 \leq t \leq \tau_1.$

Then after $\tau_1$ we restart the construction by defining the drift $(\mu_t^2)_{t \geq \tau_1}$ where $\mu_{\tau_1}^2 = \mu$ and $(\mu_t^2)_{t \geq \tau_1}$ fulfills (2.14). Set $S_{\tau_1}^2 = 0$ and

$$dS_t^2 = \mu_t^2 \, dt + dB_t, \qquad t \geq \tau_1.$$

Furthermore, we perform the Lévy transform resulting in $(M_t^2)_{t \geq \tau_1}$ and we define the stopping time

$$\tau_2 := \inf\{t > \tau_1 : |M_t^2| \geq \delta\} \wedge (\delta + \tau_1).$$

By this construction, we have that the estimate (4.3) holds for $(\mu_t^2)_{\tau_1 \leq t \leq \tau_2}$, and we may continue the construction in the same fashion.

Now we proceed formally:

Set $\tau_0 = 0$ and define recursively for $l \geq 1$ $(\tilde{S}_t^l)_{t \geq 0}$ by $\tilde{S}_0^l = 0$ and

$$d\tilde{S}_t^l = \tilde{\mu}_t^l \, dt + dW_t^l, \qquad t \geq 0,$$

where the Brownian motion $(W_t^l)_{t \geq 0}$ is given by

$$W_t^l := B_{\tau_{l-1} + t} - B_{\tau_{l-1}}, \qquad t \geq 0,$$

and the drift process $(\tilde{\mu}_t^l)_{t\geq 0}$ is given by

$$\tilde{\mu}_t^l = \frac{2\mu}{1 + \exp(2\mu W_t^l + 2\mu^2 t)} \tag{4.4}$$

[compare to (2.15)].

The integrand $(H_t^l)_{t\geq 0}$ is defined analogously to Section 3 by

$$H_t^l = \operatorname{sign}(\tilde{S}_t^l), \qquad t \geq 0,$$

and the stopping time $\gamma_l$ is defined by

$$\gamma_l := \inf\{t : |(H^l \cdot \tilde{S}^l)_t| \geq \delta\} \wedge \delta. \tag{4.5}$$

Then we set

$$\tau_l = \tau_{l-1} + \gamma_l \tag{4.6}$$

and go on with the recursive definition.

Finally, for $l \geq 1$ we introduce the processes $(\tilde{N}_t^l)_{t\geq 0}$ and $(N_t^l)_{t\geq 0}$ by

$$\tilde{N}_t^l := (H^l \cdot \tilde{S}^l)_t, \qquad t \geq 0,$$

and

$$N_t^l := \tilde{N}_{t \wedge \gamma_l}^l, \qquad t \geq 0.$$

Note that by the considerations of Sections 2 and 3 $(\tilde{N}_t^l)_{t\geq 0}$ as well as its stopped version $(N_t^l)_{t\geq 0}$ are martingales in their own filtrations for $l \geq 1$.

REMARK 4.2. It is evident that $(\gamma_l)_{l\geq 1}$ as well as $\{(W_t^l)_{0\leq t\leq \gamma_l}\}_{l\geq 1}$ are i.i.d. sequences. By $\sigma(N_t^l, t \geq 0) \subseteq \mathcal{F}_{\gamma_l}^{W^l}$, it holds that $(N_t^l)_{l\geq 1}$ are independent (and identically distributed) processes and that $\mathcal{F}_{\tau_{l-1}}^B$ is independent of $(W_t^l)_{t\geq 0}$ for $l \geq 1$. By these observations, it follows that $\mathcal{F}_{\tau_{l-1}}^B$ is independent of $(N_t^l)_{t\geq 0}$ for $l \geq 1$.

We need to show that the union of the stochastic intervals $\bigcup_{l\geq 1} [\![\tau_{l-1}, \tau_l]\!]$ equals the whole real line.

LEMMA 4.3. *Let $(\tau_l)_{l=0}^\infty$ be defined by $\tau_0 = 0$ and (4.6) for $l \geq 1$. Then*

$$\mathbb{P}[\tau_l \to \infty, \ l \to \infty] = 1.$$

PROOF. We already noticed in Remark 4.2 that the interval lengths $\tau_l - \tau_{l-1} = \gamma_l$ are positive and identically distributed. A well-known result (see, e.g., [2], Proposition 4.14) tells us that

$$\lim_{l\to\infty} \tau_l = \sum_{l=1}^\infty (\tau_l - \tau_{l-1}) = \infty$$

almost surely. □

The last step is to show that the process which is given by the concatenation of the Lévy transforms $N^l, l \geq 1$, on the respective stochastic intervals is a Brownian motion in its own filtration. We want to apply Lévy's criterion and first concentrate on proving the martingale property. We need three lemmas:

LEMMA 4.4. *Let $(G_t)_{t \geq 0}$ be a martingale in its own filtration. Then*

$$(G_{(t-x)^+})_{t \geq 0}$$

*is also a martingale in its own filtration for each fixed number $x \geq 0$.*

PROOF. Obvious. □

LEMMA 4.5. *Let $\eta$ be an $\mathcal{F}^B$-stopping time and $(G_t)_{t \geq 0}$ be a continuous martingale in its own filtration such that $\mathcal{F}^B_\eta \perp\!\!\!\perp (G_t)_{t \geq 0}$. Define the filtration $\mathcal{F}_t := \mathcal{F}^B_\eta \vee \mathcal{G}_t$, where $\mathcal{G}_t := \sigma(G_{(u-\eta)^+}, 0 \leq u \leq t)$, then $(G_{(t-\eta)^+})_{t \geq 0}$ is a martingale w.r.t. $(\mathcal{F}_t)_{t \geq 0}$.*

PROOF. We want to show that

$$\mathbb{E}[G_{(t-\eta)^+} \mid \mathcal{F}^B_\eta \vee \sigma(G_{(u-\eta)^+}, 0 \leq u \leq s)] = G_{(s-\eta)^+}.$$

Note that

$$\mathcal{F}^B_\eta \vee \sigma(G_{(u-\eta)^+}, 0 \leq u \leq s) \subseteq \mathcal{F}^B_\eta \vee \sigma(G_{(u-\cdot)^+}, 0 \leq u \leq s)$$

and consider an event in the latter sigma-algebra given by

$$A := \{B(\cdot \wedge \eta) \in C, L(\cdot) \in D\},$$

with Borel sets $C, D$ of $C[0, \infty)$ where $C[0, \infty)$ is the space of continuous functions on $[0, \infty)$ equipped with the topology of uniform convergence on compacts. We regard $B(\cdot \wedge \eta)$ as a random function

$$B(\cdot \wedge \eta) : \Omega \to C[0, \infty),$$

and

$$L(\cdot) := G_{(u_1 - \cdot)^+} : \Omega \to C[0, \infty)$$

for some $0 \leq u_1 \leq s$. Furthermore, we define the random functions

$$H_1(\cdot) := G_{(t-\cdot)^+} : \Omega \to C[0, \infty) \quad \text{and}$$
$$H_2(\cdot) := G_{(s-\cdot)^+} : \Omega \to C[0, \infty).$$

Now consider that the law $\nu$ of $(H_1, H_2, B(\cdot \wedge \eta), L, \eta)$ on the space

$$\Theta := (C[0,\infty))^4 \times (\mathbb{R}_+ \cup \{\infty\})$$

can be decomposed as

$$d\nu(x_1, x_2, x_3, x_4, x_5) = d\nu(x_1, x_2, x_4 \mid x_3, x_5) \, d\mu(x_3, x_5),$$

where $\mu$ is the law of $(B(\cdot \wedge \eta), \eta)$ and $\nu(\cdot, \cdot, \cdot \mid x_3, x_5)$ is the conditional law of $(H_1, H_2, L)$ knowing $B(\cdot \wedge \eta) = x_3$ and $\eta = x_5$. The martingale property of $G_{(t-x)^+}$ for each $x \geq 0$ (Lemma 4.4) implies

$$(4.7) \quad \int_\Theta x_1(x) 1_D(x_4(x)) \, d\vartheta(x_1, x_2, x_4) = \int_\Theta x_2(x) 1_D(x_4(x)) \, d\vartheta(x_1, x_2, x_4),$$

for each $x \geq 0$ where $\vartheta$ is the (unconditional) law of $(H_1, H_2, L)$.

Furthermore, the hypotheses of the lemma entail the independence of $H_1, H_2, L$ from $B(\cdot \wedge \eta)$ and $\eta$, so it follows that $\nu(x_1, x_2, x_4 \mid x_3, x_5)$ does not depend on $(x_3, x_5)$ and, thus,

$$(4.8) \quad d\nu(x_1, x_2, x_4 \mid x_3, x_5) = d\vartheta(x_1, x_2, x_4).$$

By the decomposition of $\nu$ and (4.8), we can write

$$\mathbb{E}[G_{(t-\eta)^+} 1_A] = \int_\Theta x_1(x_5) 1_C(x_3) 1_D(x_4(x_5)) \, d\nu(x_1, x_2, x_3, x_4, x_5)$$

$$= \int_\Theta x_1(x_5) 1_D(x_4(x_5)) \, d\vartheta(x_1, x_2, x_4) \, 1_C(x_3) \, d\mu(x_3, x_5),$$

which by (4.7) equals

$$\int_\Theta x_2(x_5) 1_D(x_4(x_5)) \, d\vartheta(x_1, x_2, x_4) \, 1_C(x_3) \, d\mu(x_3, x_5) = \mathbb{E}[G_{(s-\eta)^+} 1_A].$$

For more general sets $A$ of the form

$$A := \{B(\cdot \wedge \eta) \in C, G_{(u_1-\eta)^+} \in D_1, \ldots, G_{(u_n-\eta)^+} \in D_n\},$$

the equality

$$\mathbb{E}[G_{(t-\eta)^+} 1_A] = \mathbb{E}[G_{(s-\eta)^+} 1_A]$$

holds by the same argument, which proves the lemma. $\square$

LEMMA 4.6. *Let $X$ be a random variable, $\mathcal{G}$ and $\mathcal{H}$ be sigma-algebras in a probability space $(\Omega, \mathcal{F}, \mathbb{P})$. If $\mathcal{G} \perp\!\!\!\perp \mathcal{H}$ and $X \perp\!\!\!\perp \mathcal{H}$, then*

$$\mathbb{E}[X \mid \mathcal{G} \vee \mathcal{H}] = \mathbb{E}[X \mid \mathcal{G}].$$

PROOF. Obvious. □

Finally define the process $(S_t)_{t\geq 0}$ using $\tilde{S}^l, l \geq 1$, and $(H_t)_{t\geq 0}$ using $H^l, l \geq 1$:

$$S_t := \sum_{j=1}^{l-1} \tilde{S}^j_{\gamma_j} + \tilde{S}^l_{t-\tau_{l-1}} \qquad \text{for } \tau_{l-1} \leq t \leq \tau_l,$$

thus, the drift $(\mu_t)_{t\geq 0}$ of $(S_t)_{t\geq 0}$ is given by

(4.9) $$\mu_t := \tilde{\mu}^l_{t-\tau_{l-1}} \qquad \text{for } \tau_{l-1} \leq t \leq \tau_l,$$

and the integrand is defined as

$$H_t := H^l_{t-\tau_{l-1}} \qquad \text{for } \tau_{l-1} \leq t \leq \tau_l.$$

We obviously have

$$S_t = \int_0^t \mu_s \, ds + B_t.$$

Then the stochastic integral $(M_t)_{t\geq 0}$ is defined by

$$M_t := (H \cdot S)_t = \sum_{l=1}^{k-1} (H^l \cdot \tilde{S}^l)_{\gamma_l} + (H^k \cdot \tilde{S}^k)_{t-\tau_{k-1}} \qquad \text{for } \tau_{k-1} \leq t \leq \tau_k.$$

(4.10)
Note that, by construction, $S_t$ and $M_t$ are continuous processes.

PROPOSITION 4.7. *The process $(M_t)_{t\geq 0}$ as defined in (4.10) satisfies*

(4.11) $$M_t = \sum_{l=1}^{\infty} N^l_{(t-\tau_{l-1})^+},$$

*where the sum converges in $L^2$. $(M_t)_{t\geq 0}$ is a martingale in its own filtration. That is, for $0 \leq s < t < \infty$,*

$$\mathbb{E}[M_t | \mathcal{F}^M_s] = M_s.$$

PROOF. First we show that the sum on the right-hand side of (4.11) converges in $L^2$. Note that

$$\sum_{l=1}^{k} N^l_{(t-\tau_{l-1})^+}$$

$$= \int_0^{\tau_1 \wedge t} \text{sign}(\tilde{S}^1_{s-\tau_0})\tilde{\mu}^1_{s-\tau_0} \, ds + \cdots + \int_{\tau_{k-1} \wedge t}^{\tau_k \wedge t} \text{sign}(\tilde{S}^k_{s-\tau_{k-1}})\tilde{\mu}^k_{s-\tau_{k-1}} \, ds$$

$$+ \int_0^{\tau_1 \wedge t} \text{sign}(\tilde{S}^1_{s-\tau_0}) \, dB_s + \cdots + \int_{\tau_{k-1} \wedge t}^{\tau_k \wedge t} \text{sign}(\tilde{S}^k_{s-\tau_{k-1}}) \, dB_s$$

$$= \int_0^{\tau_k \wedge t} H_s \mu_s \, ds + \int_0^{\tau_k \wedge t} H_s \, dB_s,$$

which by Lemma 4.3 a.s. converges to

$$(4.12) \qquad M_t = \int_0^t H_s \mu_s \, ds + \int_0^t H_s \, dB_s \qquad \text{as } k \to \infty.$$

Furthermore, we know from Proposition 2.3 that

$$|H_t \mu_t| \leq 2\mu \qquad \text{for } t \geq 0.$$

By the Doob inequality and the Itô isometry, we get

$$\mathbb{E}\left[\sup_{0 \leq u \leq t}\left|\int_0^u H_s \, dB_s\right|^2\right] \leq 4\mathbb{E}\left[\left(\int_0^t H_s \, dB_s\right)^2\right] = 4t,$$

as $|H_s| = 1$. The $L^2$ convergence of the infinite sum follows.

Now we prove that $(M_t)_{t \geq 0}$ is a martingale in its own filtration. Define the filtrations $(\mathcal{G}_t^l)_{t \geq 0}$ for $l \geq 1$ by

$$\mathcal{G}_t^l := \sigma(N_{(u-\tau_{l-1})^+}^l, 0 \leq u \leq t),$$

and consider that

$$(4.13) \qquad \mathbb{E}[M_t \mid \mathcal{F}_s^M] = \mathbb{E}\left[\sum_{l=1}^\infty \mathbb{E}\left[N_{(t-\tau_{l-1})^+}^l \;\Big|\; \bigvee_{j=1}^\infty \mathcal{G}_s^j\right] \;\Big|\; \mathcal{F}_s^M\right],$$

where $L^2$-convergence allows us to exchange summation and expectation and we used that $\mathcal{F}_s^M \subseteq \bigvee_{j=1}^\infty \mathcal{G}_s^j$ for $s \geq 0$. Furthermore, notice that

$$\bigvee_{j=1}^\infty \mathcal{G}_s^j \subseteq \bigvee_{j=1}^l \mathcal{G}_s^j \vee \sigma(N^{l+1}, N^{l+2}, \ldots).$$

To see this, note that for $m \geq l+1$ we have that $\mathcal{G}_s^m \subseteq \sigma(N^m, \tau_{m-1} \wedge s)$ and that, by definition, $\sigma(\tau_{m-1} \wedge s) \subseteq \sigma(N^{m-1}, \tau_{m-2} \wedge s)$. Continuing this inductively, we get $\mathcal{G}_s^m \subseteq \sigma(N^m, \ldots, N^{l+1}, \tau_l \wedge s)$ and finally, $\sigma(\tau_l \wedge s) \subseteq \mathcal{G}_s^l$.

Thus, define $\mathcal{H}_s^l := \bigvee_{j=1}^l \mathcal{G}_s^j \vee \sigma(N^{l+1}, N^{l+2}, \ldots)$ and recall that $\bigvee_{j=1}^l \mathcal{G}_s^j \subseteq \mathcal{F}_{\tau_l \wedge s}^B$ which is independent of $N^k, k \geq l+1$; recall Remark 4.2. By Lemma 4.6 and the tower law, the inner conditional expectation in (4.13) is given by

$$\mathbb{E}\left[\mathbb{E}[N_{(t-\tau_{l-1})^+}^l \mid \mathcal{H}_s^l] \;\Big|\; \bigvee_{j=1}^\infty \mathcal{G}_s^j\right] = \mathbb{E}\left[\mathbb{E}\left[N_{(t-\tau_{l-1})^+}^l \;\Big|\; \bigvee_{j=1}^l \mathcal{G}_s^j\right] \;\Big|\; \bigvee_{j=1}^\infty \mathcal{G}_s^j\right]$$

$$= \mathbb{E}\left[N_{(t-\tau_{l-1})^+}^l \;\Big|\; \bigvee_{j=1}^l \mathcal{G}_s^j\right].$$

Note that $\bigvee_{j=1}^l \mathcal{G}_s^j \subseteq \mathcal{F}_s^l := \mathcal{F}_{\tau_{l-1}}^B \vee \mathcal{G}_s^l$. By applying Lemma 4.5 with $(G_t)_{t \geq 0} = (N_t^l)_{t \geq 0}$, we get that

$$\mathbb{E}\left[N_{(t-\tau_{l-1})^+}^l \;\Big|\; \bigvee_{j=1}^l \mathcal{G}_s^j\right] = N_{(s-\tau_{l-1})^+}^l.$$

Substituting these results into (4.13) and using the representation (4.11) for $M_s$, we get that

$$\mathbb{E}[M_t \mid \mathcal{F}_s^M] = \mathbb{E}\left[\sum_{l=1}^{\infty} N_{(s-\tau_{l-1})^+}^l \,\bigg|\, \mathcal{F}_s^M\right] = \mathbb{E}[M_s \mid \mathcal{F}_s^M] = M_s. \qquad \square$$

PROOF OF (ii) IN THEOREM 1.2. By construction and by Proposition 4.7, $(M_t)_{t\geq 0}$ is a continuous martingale and its bracket is $\langle M \rangle_t = t$ by (4.12) and by $|H_s| = 1, s \geq 0$, hence, $(M_t)_{t\geq 0}$ is a Brownian motion (in its own filtration).

For $(\mu_t)_{t\geq 0}$, the drift of $S$, we conclude that, due to (4.1), Lemma 4.1, (4.5), (4.9) and Lemma 4.3,

(4.14) $$\sup_{t\geq 0} |\mu_t - \mu| \leq \delta(3\mu^3 + 2\mu^2) \qquad \text{a.s.,}$$

which can be made arbitrarily small. $\square$

## APPENDIX: MICHEL ÉMERY'S QUESTION

We now take up again the question discussed at the end of the introduction. We adopt the notation from there but assume, slightly more generally, that the independent sequence $(\varepsilon_n)_{n\leq 0}$ of $\{-1,1\}$-valued random variables fulfill the condition

$$\mathbb{P}[\varepsilon_n = 1] = 1 - \mathbb{P}[\varepsilon_n = 1] = p_n, \qquad n \leq 0,$$

for some sequence $(p_n)_{n\leq 0}$ in $]0,1[$ satisfying

(A.1) $$\sum_{n=-\infty}^{0} \min(p_n, 1-p_n) = \infty.$$

In the sequel we call $\{-1,1\}$-valued random variables *Bernoulli* variables and let $\mathbb{N}_-$ denote the integers less than or equal to zero.

The role of the regularity condition (A.1) is explained in the following lemma.

LEMMA A.1. *The law of the $\{-1,1\}^{\mathbb{N}_-}$-valued random variable $(\varepsilon_n)_{n\leq 0}$ is diffuse iff (A.1) holds true.*

PROOF. First assume that there exists an atom $A = (a_0, a_{-1}, \ldots)$ with $\mathbb{P}[A] > 0$, then $\prod_{n=0}^{-\infty} p_n^{(1+a_n)/2}(1-p_n)^{(1-a_n)/2} > 0$ which is equivalent to

$$\sum_{n=-\infty}^{0} (1 - p_n^{(1+a_n)/2}(1-p_n)^{(1-a_n)/2}) = \sum_{n=-\infty}^{0} ((1-p_n)^{(1+a_n)/2} p_n^{(1-a_n)/2}) < \infty,$$

which implies that the sum in (A.1) is finite. On the other hand, if the sum in (A.1) is finite, this is equivalent to $\prod_{n=0}^{-\infty} p_n^{(1-a_n)/2}(1-p_n)^{(1+a_n)/2} > 0$ for a sequence $(a_n)_{n \leq 0}$ such that

$$p_n^{(1+a_n)/2}(1-p_n)^{(1-a_n)/2} = \min(p_n, 1-p_n) \qquad \text{for } n \leq 0,$$

and we find $A = (a_0, a_{-1}, \ldots)$ with $\mathbb{P}[A] > 0$. $\square$

We call a $\{-1, 1\}$-valued random variable $X$ *symmetric Bernoulli* if

$$\mathbb{P}[X = 1] = \mathbb{P}[X = -1] = \tfrac{1}{2}.$$

LEMMA A.2. *Let $(\varepsilon_n)_{n \leq 0}$ be a sequence of Bernoulli random variables, and $(h_n)_{n \leq 0}$ an i.i.d. sequence of symmetric Bernoulli variables independent of $(\varepsilon_n)_{n \leq 0}$. Then:*

(a) $(h_n \varepsilon_n)_{n \leq 0}$ *is an i.i.d. sequence of symmetric Bernoulli random variables and*

(b)
$$\mathrm{law}[(\varepsilon_n)_{n \leq 0} | (h_n \varepsilon_n)_{n \leq 0}] = \mathrm{law}[(\varepsilon_n)_{n \leq 0}] \qquad a.s.$$

PROOF. Fix $N \geq 1$ and consider signs $x_1, \ldots, x_N$ as well as indices $i_1, \ldots, i_N$. Then by independence of $(h_n)_{n \leq 0}$ and $(\varepsilon_n)_{n \leq 0}$ combined with the i.i.d. property, we get

$$\mathbb{P}[h_{i_1} \varepsilon_{i_1} = x_1, \ldots, h_{i_N} \varepsilon_{i_N} = x_N]$$
$$= \sum_{y_1, \ldots, y_n} \mathbb{P}[h_{i_1} = x_1/y_1, \ldots, h_{i_N} = x_N/y_N] \mathbb{P}[\varepsilon_{i_1} = y_1, \ldots, \varepsilon_{i_N} = y_N]$$
$$= \left(\frac{1}{2}\right)^N \sum_{y_1, \ldots, y_n} \mathbb{P}[\varepsilon_{i_1} = y_1, \ldots, \varepsilon_{i_N} = y_N] = \left(\frac{1}{2}\right)^N,$$

which proves (a). For proving (b), we fix again $N, M \geq 1$ and consider signs $x_1, \ldots, x_N$ and $y_1, \ldots, y_N$ as well as indices $i_1, \ldots, i_N$ such that $\mathbb{P}[h_{i_1} \varepsilon_{i_1} = x_1, \ldots, h_{i_N} \varepsilon_{i_N} = x_N] > 0$. By independence of $(h_n)_{n \leq 0}$ and $(\varepsilon_n)_{n \leq 0}$ and the previous argument, we can calculate that

$$\mathbb{P}[\varepsilon_{i_1} = y_1, \ldots, \varepsilon_{i_N} = y_N \mid h_{i_1} \varepsilon_{i_1} = x_1, \ldots, h_{i_N} \varepsilon_{i_N} = x_N]$$
$$= \frac{\mathbb{P}[\varepsilon_{i_1} = y_1, \ldots, \varepsilon_{i_N} = y_N, h_{i_1} \varepsilon_{i_1} = x_1, \ldots, h_{i_N} \varepsilon_{i_N} = x_N]}{\mathbb{P}[h_{i_1} \varepsilon_{i_1} = x_1, \ldots, h_{i_N} \varepsilon_{i_N} = x_N]}$$
$$= (1/2)^N \frac{\mathbb{P}[\varepsilon_{i_1} = y_1, \ldots, \varepsilon_{i_N} = y_N]}{(1/2)^N}$$
$$= \mathbb{P}[\varepsilon_{i_1} = y_1, \ldots, \varepsilon_{i_N} = y_N],$$

which proves (b). □

Assuming (A.1), we can find disjoint, infinite subsets $(I_n)_{n\leq 0}$ of $\mathbb{N}_-$ such that $i > n$, for all $i \in I_n$, and

(A.2) $$\sum_{i \in I_n} \min(p_i, 1 - p_i) = \infty.$$

For these sets we define the following infinite sequence $(I^{(l)})_{l=0}^\infty$ of subsets of $\mathbb{N}_-$ by

(A.3)
$$I^{(0)} := I_0,$$
$$I^{(1)} = \bigcup_{n \in I^{(0)}} I_n,$$
$$I^{(2)} = \bigcup_{n \in I^{(1)}} I_n,$$
$$\vdots \quad \text{etc.}$$

Additionally to the sequence $(I^{(l)})_{l=0}^\infty$ from (A.3), we furthermore introduce

$$J = \mathbb{N}_- \setminus \left( \{0\} \cup \bigcup_{l=0}^\infty I^{(l)} \right).$$

In the following lemma we summarize three properties of these sets.

LEMMA A.3. *For the sequence $(I^{(l)})_{l=0}^\infty$ defined in (A.3) we have the following:*

(a) $I^{(l)} \subset \{\ldots, -l-2, -l-1\}$ *for* $l \geq 0$.
(b) *The sets* $(I^{(l)})_{l=0}^\infty$ *are mutually disjoint.*
(c) *For each* $m \in J$ *we have* $I_m \subset J$.

PROOF. *Proof of (a)*: We prove the statement by induction. $0 \notin I^{(0)}$ by construction. Thus, assume that the statement holds for $I^{(0)}, \ldots, I^{(n)}$. For $I^{(n+1)}$ consider that

$$I^{(n+1)} = \bigcup_{x \in I^{(n)}} I_x,$$

and by the induction hypothesis $x \leq -n - 1$ for $x \in I^{(n)}$. But then also $y \leq x - 1$, for all $y \in I_x$, thus, it follows that $y \leq -n - 2$ and $I_x \subset \{\ldots, -n - 3, -n - 2\}$ for each $x \in I^{(n)}$, which proves (a).

*Proof of (b)*: Again by induction, let us assume that $I^{(0)}, \ldots, I^{(n)}$ are pairwise disjoint. We want to prove that $I^{(0)}, \ldots, I^{(n+1)}$ are also pairwise disjoint. Take $m \in I^{(n+1)} = \bigcup_{x \in I^{(n)}} I_x$. If we had $m \in I^{(j)}$ for some $1 \le j \le n$, then $m \in I_y$ for some $y \in I^{(j-1)}$ and also $m \in I_w$ for some $w \in I^{(n)}$. But as the $I_i, i \in \mathbb{N}_-$ are disjoint, this implies $y = w$ so that $I^{(n)} \cap I^{(j-1)} \ne \varnothing$, which is a contradiction. Finally, if $m \in I^{(0)}$, then $w = 0$, but $0 \notin I^{(i)}$ for $i \ge 1$ by (a).

*Proof of (c)*: Let $m \in J$ and $x \in I_m$. Let $I^{(-1)} := \{0\}$. Assume that there is a $k \ge 0$ such that $x \in I^{(k)}$; this implies that there is $y \in I^{(k-1)}$ such that $x \in I_y$. Then $m = y$ and $m \in I^{(k-1)}$, which is a contradiction. $\square$

We now can prove a positive answer for M. Émery's question.

THEOREM A.4. *Let $(\varepsilon_n)_{n \le 0}$ be a sequence of independent $\{-1, 1\}$-valued random variables such that*

(A.4) $$\sum_{n=-\infty}^{0} \min(\mathbb{P}[\varepsilon_n = 1], \mathbb{P}[\varepsilon_n = -1]) = \infty.$$

*Then there is a predictable process $(h_n)_{n \le 0}$ of $\{-1, 1\}$-valued random variables, such that $(h_n \varepsilon_n)_{n \le 0}$ is an i.i.d. sequence of symmetric Bernoulli random variables.*

PROOF. Consider the disjoint, infinite subsets $(I_n)_{n \le 0}$ of $\mathbb{N}_-$ verifying (A.2).

By Lemma A.1, we may find Borel-functions

$$f_n : \{-1, 1\}^{\mathbb{N}_-} \to \{-1, 1\},$$

such that

(A.5) $$h_n = f_n((\varepsilon_i)_{i \in I_n})$$

satisfies

$$\mathbb{P}[h_n = 1] = \mathbb{P}[h_n = -1] = \tfrac{1}{2}.$$

We claim that these $(h_n)_{n \le 0}$ do the job.

For this aim we show that

(A.6) $\quad \mathbb{P}[h_0 \varepsilon_0 = 1 | (h_n \varepsilon_n)_{n \le -1}] \stackrel{\text{a.s.}}{=} \mathbb{P}[h_0 \varepsilon_0 = -1 | (h_n \varepsilon_n)_{n \le -1}] \stackrel{\text{a.s.}}{=} 1/2.$

To see the dependence structure of the $(h_n)_{n \in \mathbb{N}_-}$, note that by (A.5) for $n \ge 0$ the $\{-1, 1\}^{\mathbb{N}_-}$-valued random variable $(h_i)_{i \in I^{(n)}}$ does depend on $(\varepsilon_i)_{i \in I^{(n+1)}}$ but it is independent of $(\varepsilon_i)_{i \in I^{(n)}}$ as $I^{(n)} \cap I^{(n+1)} = \varnothing$, by Lemma A.3(b).

The $(h_i)_{i \in I_0}$ are an i.i.d. sequence of symmetric Bernoulli random variables independent of $(\varepsilon_i)_{i \in I_0}$ and $(\varepsilon_i)_{i \in J}$, hence, by Lemma A.2(b),

(A.7) $\qquad \text{law}[(\varepsilon_i)_{i \in I_0} | (h_i \varepsilon_i)_{i \in I_0}, (\varepsilon_i)_{i \in J}] \stackrel{\text{a.s.}}{=} \text{law}[(\varepsilon_i)_{i \in I_0}].$

It follows that

$$\mathbb{P}[h_0 = 1 \mid (h_i\varepsilon_i)_{i \in I_0}, (\varepsilon_i)_{i \in J}] \stackrel{\text{a.s.}}{=} \mathbb{P}[f_0((\varepsilon_i)_{i \in I_0}) = 1 \mid (h_i\varepsilon_i)_{i \in I_0}, (\varepsilon_i)_{i \in J}]$$
$$\stackrel{\text{a.s.}}{=} \mathbb{P}[f_0((\varepsilon_i)_{i \in I_0}) = 1] = \tfrac{1}{2}. \quad (A.8)$$

Similarly, we get for $(h_i)_{i \in I_0} = (f_i(\varepsilon_n, n \in I_i))_{i \in I_0}$ that

$$(A.9) \quad \text{law}[(h_i)_{i \in I_0} \mid (h_i\varepsilon_i)_{i \in I^{(1)}}, (\varepsilon_i)_{i \in J}] \text{ is a.s. i.i.d. symmetric Bernoulli.}$$

Now we claim that

$$(A.10) \quad \text{law}[h_0 \mid (h_i\varepsilon_i)_{i \in I_0 \cup I^{(1)}}, (\varepsilon_i)_{i \in J}] \text{ is a.s. symmetric Bernoulli.}$$

Indeed, fix finite index sets $K \subset I_0$, $L \subset I^{(1)}$, $M \subset J$ and signs $(x_i)_{i \in K}$ as well as $(y_i)_{i \in L}$, $(z_i)_{i \in M}$. Then

$$\mathbb{P}[h_0 = 1 \mid (h_i\varepsilon_i)_{i \in K} = (x_i)_{i \in K}, (h_i\varepsilon_i)_{i \in L} = (y_i)_{i \in L}, (\varepsilon_i)_{i \in M} = (z_i)_{i \in M}]$$
$$= \frac{\mathbb{P}[h_0 = 1, (h_i\varepsilon_i)_{i \in K} = (x_i)_{i \in K} \mid (h_i\varepsilon_i)_{i \in L} = (y_i)_{i \in L}, (\varepsilon_i)_{i \in M} = (z_i)_{i \in M}]}{\mathbb{P}[(h_i\varepsilon_i)_{i \in K} = (x_i)_{i \in K} \mid (h_i\varepsilon_i)_{i \in L} = (y_i)_{i \in L}, (\varepsilon_i)_{i \in M} = (z_i)_{i \in M}]}.$$

For the denominator consider that by (A.9) together with independence of $(\varepsilon_i)_{i \in I_0}$ from $(\varepsilon_i)_{i \in I^{(1)} \cup I^{(2)} \cup J}$ and by Lemma A.2(a),

$$\mathbb{P}[(h_i\varepsilon_i)_{i \in K} = (x_i)_{i \in K} \mid (h_i\varepsilon_i)_{i \in L} = (y_i)_{i \in L}, (\varepsilon_i)_{i \in M} = (z_i)_{i \in M}] = 2^{-|K|}.$$

By (A.8) and (A.9), we get for the numerator

$$\mathbb{P}[h_0 = 1, (h_i\varepsilon_i)_{i \in K} = (x_i)_{i \in K} \mid (h_i\varepsilon_i)_{i \in L} = (y_i)_{i \in L}, (\varepsilon_i)_{i \in M} = (z_i)_{i \in M}]$$
$$= \mathbb{P}[h_0 = 1, (h_i\varepsilon_i)_{i \in K} = (x_i)_{i \in K}] = \tfrac{1}{2} 2^{-|K|},$$

and (A.10) follows as the same conclusion passes to infinite index sets.

Continuing analogously, we get by induction that

$$\text{law}[h_0 \mid (h_i\varepsilon_i)_{i \in \bigcup_{l=0}^{\infty} I^{(l)}}, (\varepsilon_i)_{i \in J}] \text{ is a.s. symmetric Bernoulli,}$$

which gives the claim (A.6) since $\sigma(h_n\varepsilon_n, n \leq -1)$ is contained in $\sigma(\varepsilon_n h_n, n \in \bigcup_{l=0}^{\infty} I^{(l)}, \varepsilon_n, n \in J)$ by Lemma A.3(c).

Analogous arguments show that

$$\mathbb{P}[h_i\varepsilon_i = 1 \mid (h_n\varepsilon_n)_{n \leq i-1}] \stackrel{\text{a.s.}}{=} \mathbb{P}[h_i\varepsilon_i = -1 \mid (h_n\varepsilon_n)_{n \leq i-1}] \stackrel{\text{a.s.}}{=} 1/2,$$

for any $i \leq -1$, which proves the theorem. $\square$

**Acknowledgments.** We want to thank Marc Yor for introducing us to this problem and for repeated discussions and Michel Émery for posing and extensively discussing the question in discrete time. Special thanks go to Josef Teichmann for fruitful discussions, to Vilmos Prokaj for realizing the explicit form of $\mu_t$ and to Pavel Gapeev for bringing the reference [4] to our attention. Furthermore, we wish to thank Uwe Schmock for pointing out an inaccuracy in a draft of this article and an anonymous referee for a very helpful report.

COMPUTER AND AUTOMATION RESEARCH
 INSTITUTE OF THE HUNGARIAN
 ACADEMY OF SCIENCES
KENDE UTCA 13-17
1111 BUDAPEST
HUNGARY
E-MAIL: rasonyi@sztaki.hu

FACULTY OF MATHEMATICS
UNIVERSITY OF VIENNA
NORDBERGSTRASSE 15
1090 WIEN
AUSTRIA
E-MAIL: walter.schachermayer@univie.ac.at

RAIFFEISEN CAPITAL MANAGEMENT
SCHWARZENBERGPLATZ 3
1010 WIEN
AUSTRIA
E-MAIL: richard.warnung@rcm.at